\documentclass{article}
\usepackage{amsfonts}

\input{tcilatex}

\begin{document}

\title{CLASSE D'EQUIVALENCE FORMELLE D'UN D-MODULE D'AIRY}
\author{LOTFI\ SA\"{I}DANE \\
%EndAName
D\'{e}partement de Math\'{e}matiques \\
Facult\'{e} des Sciences de Tunis. \\
2092, EL-MANAR\ I, TUNIS \\
TUNISIE }
\maketitle

\textbf{R\'{e}sum\'{e} : Dans ce travail on \'{e}tudie une classe d'\'{e}%
quivalence formelle d'un D-module d'Airy, qu'on compare avec une classe d'%
\'{e}quivalence analytique. On calcule les facteurs d\'{e}terminants d'un op%
\'{e}rateur d'Airy de bidegr\'{e} }$\mathbf{(n,m)}$\textbf{\ et on pr\'{e}%
cise les coefficients de l'op\'{e}rateur qui interviennent dans la d\'{e}%
termination des facteurs d\'{e}terminants. On donne enfin, le mod\`{e}le
canonique, voir [BV1], d'un op\'{e}rateur d'Airy de bidegr\'{e} }$\mathbf{%
(n,m)}$\textbf{.}

\bigskip

Soit $k$ un corps de caract\'{e}ristique z\'{e}ro alg\'{e}briquement clos, $%
\mathbf{A}_{k}^{1}$= $Spec(k[x])$ la ligne affine sur $k$, $\partial=\frac {d%
}{dx}$, $\mathcal{D}$ $=$ $k$[$x$,$\partial$] l'alg\`{e}bre de Weyl, si $%
P(x)=$ $\Sigma_{i=0}^{n}a_{i}x^{i}$ $\in$ $k[x]$, alors $P$($\partial
)=\Sigma_{i=0}^{n}$ $a_{i}\partial^{i}$ est l'op\'{e}rateur diff\'{e}rentiel
correspondant (\`{a} coefficients constants) sur $\mathbf{A}_{k}^{1}$\textbf{%
.}

On appelle op\'{e}rateur d'Airy de bidegr\'{e} $(n,m)$ tout op\'{e}rateur
diff\'{e}rentiel de $\mathcal{D}$ de la forme $P_{n}(\partial)$ $+$ $Q_{m}(x)
$, o\`{u} $P_{n}(x)$ et $Q_{m}(x)$ sont deux polyn\^{o}mes de $k[x]$ de degr%
\'{e}s respectifs $n$ et $m$. Le $D$-module correspondant $\mathcal{D} $/$%
\mathcal{D}$($P_{n}(\partial)$ $+Q_{m}(x))$ est appel\'{e} $D$-module d'Airy
de bidegr\'{e} $(n,m)$.

Dans la suite on prendra $k$ = $\mathbb{C}$ et $K$ = $\mathbb{C}((x))$ le
corps des s\'{e}ries formelles \`{a} coefficients complexes. Si $V_{1}$et $%
V_{2}$ sont deux $D$-modules d'Airy de bidegr\'{e} $(n,m)$, $M_{1}$ et $%
M_{2} $ les syst\'{e}mes diff\'{e}rentiels lin\'{e}aires d'ordre $1$ associ%
\'{e}s; on dira que $V_{1}$ est formellement \'{e}quivalent \`{a} $V_{2}$ si
et seulement si $\exists\mathrm{\ }P\in$ $Gl(n,$ $K)$ telle que :

\begin{center}
$M_{2}$ = $P^{-1}$ $M_{1}P$ $-P^{-1}\partial P.$
\end{center}

\noindent Ce qui veut dire que si $Y$ est une solution formelle de $\partial
Y$ $=M_{2}Y$ alors $Z=PY$ est une solution formelle de $\partial Z=M_{1}Z$.

Une classe d'\'{e}quivalence formelle contient une classe d'\'{e}quivalence
analytique.

On se propose d'\'{e}tudier une classe d'\'{e}quivalence formelle d'un $D$%
-module d'Airy. Dans le paragraphe 1 on donne l'expression des facteurs d%
\'{e}terminants d'un op\'{e}rateur d'Airy de bidegr\'{e} $(n,m)$ et les
coefficients de l'op\'{e}rateur qui d\'{e}terminent ces facteurs, on calcule 
\'{e}galement les valeurs propres de la monodromie formelle. Dans le
paragraphe 2 on donne le mod\`{e}le canonique [BV1] d'un op\'{e}rateur
d'Airy.\medskip

\textbf{1. }\underline{\textbf{Facteurs d\'{e}terminants}}.\medskip

\textbf{1.1. }

Soit $K=%
%TCIMACRO{\U{2102} }%
%BeginExpansion
\mathbb{C}
%EndExpansion
((x))$ le corps des s\'{e}ries formelles \`{a} coefficients complexes muni
de sa d\'{e}rivation usuelle $\partial$ = $\frac{d}{dx}$ et de sa valuation $%
x$-adique $v$. Soit $\overline{K}$ la cl\^{o}ture alg\'{e}brique de $K$, la
valuation et la d\'{e}rivation s'\'{e}tendent de fa\c{c}on unique \`{a} $%
\overline{K}.$

Soit

\[
L\mathbf{=\partial}^{n}+\sum_{i=0}^{n}a_{i}\partial^{i}\qquad(a_{i}\in K) 
\]

\noindent un op\'{e}rateur diff\'{e}rentiel \`{a} coefficients dans $K$. On
suppose que les $a_{i}$ sont analytiques dans un voisinage de z\'{e}ro sauf
en z\'{e}ro et que z\'{e}ro est un p\^{o}le d'ordre fini des $a_{i}$, le th%
\'{e}or\`{e}me de Hukuhara et Turrittin montre l'existence d'une base de $n$
solutions, de l'\'{e}quation diff\'{e}rentielle $L(u)=0$, $u_{1}$,...,$u_{n}$
de la forme:

\[
u_{i}(x)=\exp Q_{i}(x)x^{\lambda_{i}}v_{i}(x)\qquad(i=1,...,\text{ }n) 
\]

\noindent o\`{u} $Q_{i}(x)$ est un polyn\^{o}me en $x^{-1/e}$ sans terme
constant ($e$ entier naturel), $\lambda_{i}$ est un nombre complexe et:

\[
v_{i}(x)=\sum_{j=0}^{n_{i}}v_{i,j}(x)(Logx)^{j} 
\]

\noindent o\`{u}

\[
v_{i,j}=\sum_{k=0}^{\infty}v_{i,j,k}x^{k/e}\quad\in\overline{K} 
\]

\noindent Les $Q_{i}$ sont appel\'{e}s facteurs d\'{e}terminants de $L$, $%
exp(2i\pi\lambda_{j})$ est une valeur propre de la partie de la monodromie
formelle de $L$ attach\'{e}e \`{a} $Q_{j}$. La multiplicit\'{e} d'un facteur
d\'{e}terminant $Q,$ est par d\'{e}finition, le nombre d'indice $i$ v\'{e}%
rifiant $Q_{i}=Q.$ Si, il existe $j\in\{1,$..., $n\}$ v\'{e}rifiant: $%
v(a_{i})-i\geq$ $v(a_{j})-j,$ pour tout $i$ dans $\{1,$..., $n\},$ et $%
v(a_{i})-i>$ $v(a_{j})-j,$ pour $i>j.$ Alors, $n-j$ est appell\'{e} "indice
caract\'{e}ristique" de $L$ et $j$ repr\'{e}sente la multiplicit\'{e} de $0$
comme facteur d\'{e}terminant de $L.$ On dit qu'un facteur d\'{e}terminant $%
Q $ $\neq$~$0$ est de caract\'{e}ristique $k$-uple pour $L$, si, il existe
exactement $k$ facteurs d\'{e}terminants $Q_{i_{1}},$ ..., $Q_{i_{k}}$ de $L$
v\'{e}rifiant:

\begin{center}
$v(Q_{i_{j}}-Q)$ $>$ $v(Q)$\qquad$(j=1,...,k)$
\end{center}

\noindent L'op\'{e}rateur $L$ est dit de caract\'{e}ristique $k$-uple s'il
existe un facteur d\'{e}terminant de caract\'{e}ristique $k$-uple et si les
autres sont de caract\'{e}ristique $s$-uple avec $s$ $\leq$ $k$. Autrement,
soit:

\begin{center}
$P_{L}(x,X)$ $=$ $\Sigma_{j=0}^{n}a_{j}(x)X^{j}$ \quad($a_{n}$ $=1$)
\end{center}

\noindent Le symbole de $L$. Il est clair qu'il existe $n$ fonctions
complexes, $\xi_{1}(x)$,...,\textrm{\ }$\xi_{n}(x),$dans $\bar{K},$ continue
au voisinage de z\'{e}ro sauf en z\'{e}ro tels que:

\begin{center}
$P_{L}(x,X)$ $=\dprod \limits_{j=1}^{n}(X-\xi_{j}(x))$
\end{center}

\noindent La caract\'{e}ristique de $L$ est simple si on a le r\'{e}sultat
suivant: Soit $\xi_{j}(x)$ et $\xi_{k}(x)$ deux branches non born\'{e}es de z%
\'{e}ro de $P_{L}$ tels que $x\xi_{j}(x)$ est born\'{e}e et $%
\xi_{j}(x)/\xi_{k}(x)$ $\longrightarrow1$ quand $x\longrightarrow0$, alors $%
j=k$.

Soit $\xi_{j}(x)=$ $\sum_{k=N(j)}^{\infty}\alpha_{j,k}x^{k/q}$ $(j=1,...,n)$
le d\'{e}veloppement formelle de Puiseux de $\xi_{j}(x)$, $%
\alpha_{j,N(j)}\neq$ $0$, soit $I$ $={\{}j$ tel que $N(j)/q<-1{\}}$, la
caract\'{e}ristique de $L$ est simple si les $\alpha_{j,N(j)}$, pour $j\in I,
$ sont deux \`{a} deux distincts, un r\'{e}sultat classique (voir Forsyth
[F], Helfer et Kannai [H,K]) donne les facteurs d\'{e}terminants suivant la
formule:

\begin{center}
$\frac{dQ_{j}}{dx}=\sum_{k=N(j)}^{-q-1}\alpha_{j,k}$ $x^{k/q}$
\end{center}

\noindent(Si $\xi_{j}(x)=O(1/x)$ quand $x\longrightarrow0$ alors $Q_{j}(x)=$ 
$0$, ceci est en g\'{e}n\'{e}ral vraie m\^{e}me si l'hypoth\`{e}se de la
caract\'{e}ristique n'est pas v\'{e}rifi\'{e}e).

La caract\'{e}ristique de $L$ est au plus double si on a le r\'{e}sultat
suivant:

\noindent Soient $\xi_{j}(x),$ $\xi_{k}(x)$ et $\xi_{h}(x)$ trois branches
non born\'{e}es de z\'{e}ro de $P_{L},$ tels que $x$ $\xi_{j}(x)$ est non
born\'{e}e quand $x\longrightarrow0$ et $\lim_{x\longrightarrow0}\xi
_{j}(x)/\xi_{k}(x)$ $=\lim_{x\longrightarrow0}$ $\xi_{j}(x)/\xi_{h}(x)=1,$
alors ou bien $j=k,$ ou bien $j=h,$ ou bien $k=h$.

\noindent Dans ces conditions Helfer et Kannai [H, K] ont donn\'{e}
l'expression des facteurs d\'{e}terminants de $L,$ selon le proc\'{e}d\'{e}
suivant. Soit $\tilde{P}_{L}(x,X)$ d\'{e}fini par:

\begin{center}
$\tilde{P}_{L}(x,X)$ $=P_{L}(x,X)$ $-\frac{1}{2}\frac{\partial^{2}}{\partial
x\partial X}P_{L}(x,X)$
\end{center}

\noindent On d\'{e}signe par $n-r$ l'indice caract\'{e}ristique de $L.$
Alors la d\'{e}riv\'{e}e des facteurs d\'{e}terminants non nuls (moyennant
une permutation), $Q_{r+1},$ ..., $Q_{n}$ de $L,$ sont donn\'{e}s par la
formule:

\begin{center}
$\frac{dQ_{j}}{dx}=\sum_{k=M(j)}^{-q-1}\beta_{j,k}x^{k/q},$ $j=r+1,$ ..., $n,
$
\end{center}

\noindent o\`{u}

\begin{center}
$\eta_{j}(x)$ $=$ $\sum_{k=M(j)}^{\infty}\beta_{j,k}$ $x^{k/q}$ $(j=1,...,n)$
\end{center}

\noindent est le d\'{e}veloppement formelle de Puiseux d'une branche $\eta
_{j}(x)$ de z\'{e}ro de $\tilde{P}_{L}(x,X).$

\noindent Helfer et Kannai [H,K] ont donn\'{e} l'expression des facteurs d%
\'{e}terminants d'un op\'{e}rateur diff\'{e}rentiel de caract\'{e}ristique
au plus double. Yebbou [Y] a repris le calcul des facteurs d\'{e}terminants
d'un op\'{e}rateur d'ordre $\leq$~$3$.\medskip

\textbf{1.2. }\underline{\textbf{Op\'{e}rateur d'Airy}}\textbf{:\medskip}

Soient $P_{n}$ et $Q_{m}$ deux polyn\^{o}mes de $%
%TCIMACRO{\U{2102} }%
%BeginExpansion
\mathbb{C}
%EndExpansion
\lbrack x]$ de degr\'{e}s respectifs $n$ et $m$. On pose:

\[
P_{n}=\sum_{i=1}^{n}a_{i}x^{i}\quad(a_{n}=1),\quad
Q_{m}=\sum_{j=0}^{m}b_{j}x^{j}, 
\]

\[
\partial=\frac{d}{dx},\quad\mathcal{L}(\partial,x)=\sum_{i=1}^{n}a_{i}%
\partial^{i}-\sum_{j=0}^{m}b_{j}x^{j}. 
\]

\noindent On d\'{e}signe par $D=z\frac{d}{dz}$, avec $z=\frac{1}{x}.$

\noindent On a: $D=-x$ $\partial.$

Comme, pour $k\in%
%TCIMACRO{\U{2115} }%
%BeginExpansion
\mathbb{N}
%EndExpansion
_{>0},$ $x^{k}(\frac{d}{dx})^{k}$ $=$ $x\partial(x\partial-1)$ ... $%
(x\partial-k+1),$ on a:

\begin{center}
$x^{n}\mathcal{L}(\partial,x)=$ $\sum_{i=1}^{n}a_{i}x^{n}\partial^{i}-%
\sum_{i=0}^{m}b_{i}$ $x^{i+n}.$
\end{center}

\noindent La forme de Frobenius Fuchs de $(-1)^{n}$ $x^{n}$ $\mathcal{L}%
(\partial,x)$ est:

\begin{center}
$L(D,z)=$ $\sum_{k=0}^{n}c_{k}D^{n-k},$
\end{center}

\noindent avec:

\[
c_{n}=\frac{(-1)^{n-1}}{z^{n}}Q_{m}(\frac{1}{z}) 
\]

\[
c_{k}=\sum_{i=0}^{k}(-1)^{i}\frac{a_{n-i}}{z^{i}}\sigma_{k-i,n-i}\quad\text{%
pour }k=0,...,n-1, 
\]

\noindent o\`{u}:

$\sigma_{o,j}$ $=1$, $\sigma_{h,j}=\sum_{1\leq r_{1}<r_{2}<...<r_{h}\leq
j-1}r_{1}r_{2}...r_{h}.$

\noindent On d\'{e}signe par $P_{L},$ le symbole de $L,$on a:

\begin{center}
$P_{L}=$ $\sum_{k=0}^{n}$ $c_{k}$ $X^{n-k}$
\end{center}

\noindent Le polyg\^{o}ne de Newton de $P_{L}$ admet une seule pente \`{a}
l'infini $\mu=\frac{n+m}{n}$. On d\'{e}duit, alors, que si $\xi$ est une
branche de z\'{e}ro de $P_{L}$, son d\'{e}veloppement de Puiseux est donn%
\'{e} par:

\[
\xi=\sum_{k=0}^{\infty}\frac{\alpha_{k}}{z^{1+\frac{m-k}{n}}},\quad\alpha
_{0}\neq0 
\]
\medskip

\textbf{1.2.1.}\underline{\textbf{\ Proposition}}: $L(D,z)$ est de caract%
\'{e}ristique simple et ses facteurs d\'{e}terminants sont donn\'{e}s par:

\[
Q(z)=\frac{-1}{z^{1+\frac{m}{n}}}\sum_{k=0}^{m+n-1}\frac{n\alpha_{k}}{n+m-k}%
z^{\frac{k}{n}} 
\]
\medskip

\underline{\textbf{Preuve}}: Soit $L^{1}(\frac{d}{dz},$ $z)$ $=$ $\frac {1}{%
z^{n}}L(D,z)$ $=$ $\sum_{k=0}^{n}h_{k}(z)(\frac{d}{dz})^{n-k}$. Les facteurs
d\'{e}terminants de $L$ et $L^{1}$ sont les m\^{e}mes, l'indice caract\'{e}%
ristique aussi. Soit $P_{L^{1}}(z,X)$ et $P_{L}(z,X)$ les symboles
respectifs de $L^{1}$ et $L$, $\eta(z)$ une branche de z\'{e}ro de $P_{L^{1}}
$, et $\xi(z)$ une branche de z\'{e}ro de $P_{L}$, on note $\xi_{<j}$ la
somme des mon\^{o}mes de valuation strictement inf\'{e}rieur \`{a} $j$ dans
le d\'{e}veloppement formelle de Puiseux de $\xi.$ On a (moyennant une
permutation),

\begin{center}
$\xi_{<0}$ $=z(\eta_{<-1})$
\end{center}

\noindent en comparant les termes de plus petite valuation dans $P_{L}(z$,$%
\xi(z))$ $=0,$ on a:

\begin{center}
$\alpha_{0}^{n}=(-1)^{n}$ $b_{m}$
\end{center}

\noindent L'op\'{e}rateur $L$ est donc de caract\'{e}ristique simple et les
facteurs d\'{e}terminants sont donn\'{e}s par:

\begin{center}
$z\frac{dQ(z)}{dz}=\xi_{<0}$ $=$ $\frac{1}{z^{1+\frac{m}{n}}}\sum
_{k=0}^{m+n-1}\alpha_{k}$ $z^{k/n.}$
\end{center}

\noindent D'o\`{u} r\'{e}sultat.\medskip

Soit $\xi$ une branche de z\'{e}ro de $P_{L}$, le d\'{e}veloppement de
Puiseux de $\xi$ est donn\'{e} par:

\begin{center}
$\xi(z)=$ $\sum_{k=0}^{\infty}\frac{\alpha_{k}}{z^{1+\frac{m-k}{n}}}=$ $%
\frac{1}{z^{1+\frac{m}{n}}}\sum_{k=0}^{\infty}\alpha_{k}z^{\frac{k}{n}}$

$P_{L}(z,$ $\xi(z))$ $=\sum_{k=0}^{n}$ $c_{k}\xi(z)^{n-k}$ $=0$
\end{center}

\noindent On a:

\begin{center}
$\xi(z)^{n-k}=\frac{1}{z^{n-k+(n-k)\frac{m}{n}}}\sum_{j=0}^{\infty}%
\beta_{j,n-k}$ $z^{j/n}$
\end{center}

\noindent avec

\[
\begin{array}{cc}
\beta_{j,1}= & \alpha_{j} \\ 
\beta_{0,k}= & \alpha_{0}^{k} \\ 
\beta_{j,k+1}= & \sum_{s=0}^{j}\alpha_{s}\beta_{j-s,k}%
\end{array}
\]

\noindent en comparant les termes de plus petite valuation dans $%
P_{L}(z,\xi(z))=0$ on a:

\begin{center}
$v(c_{k}\xi^{n-k})$ $=-n-m+\frac{km}{n}$

\noindent$v(c_{n})$ $=-n-m$.
\end{center}

\noindent\underline{\textbf{1}}\noindent$^{er}$\underline{\textbf{\ cas}}: $%
m=qn$ $(q$ $\geq$ $1)$.

Si $0<k<m+n$ et $k$ $\notin n%
%TCIMACRO{\U{2115} }%
%BeginExpansion
\mathbb{N}
%EndExpansion
,$ alors $\alpha_{k}$ $=0$, en effet $\beta_{k,n}$ est le coefficient de $%
\frac{1}{z^{m+n-\frac{k}{n}}}$ dans $P_{L}(z,$ $\xi(z))$ d'ou $\beta_{k,n}$ $%
=0$.

\noindent De m\^{e}me:

\begin{center}
$\{%
\begin{array}{cc}
\beta_{sn,n}= & (-1)^{n}b_{m-s}\quad si\quad0\leq s\leq q-1 \\ 
\beta_{qn,n}= & \alpha_{n-1}\beta_{0,n-1}+(-1)^{n}b_{m-q}%
\end{array}
$
\end{center}

\noindent Soit ($\alpha_{o}$, $\alpha_{n}$,..., $\alpha_{qn})$ une solution
du syst\`{e}me :

\begin{center}
$(\mathcal{A)}$ $\{%
\begin{array}{cc}
\beta_{sn,n}= & (-1)^{n}b_{m-s}\quad si\quad0\leq s\leq q-1 \\ 
\beta_{qn,n}= & \alpha_{n-1}\beta_{0,n-1}+(-1)^{n}b_{m-q}%
\end{array}
$
\end{center}

\noindent on a:

\begin{center}
$\alpha_{0}^{n}=(-1)^{n}b_{m}$, $\alpha_{n}$ $=(-1)^{n}\frac{b_{m-1}}{%
n\alpha_{0}^{n-1}},$

$\alpha_{2n}=(-1)^{n}\frac{b_{m-2}}{n\alpha_{0}^{n-1}}-\frac{n-1}{2n^{2}}%
\frac{b_{m-1}^{2}}{\alpha_{0}^{2n-1}},$ etc...\medskip
\end{center}

\textbf{1.2.2. }\underline{\textbf{Proposition}} :

Si $m=qn$ $(q\geq1)$, les facteurs d\'{e}terminants sont donn\'{e}s par:

\begin{center}
$Q(z)=$ $\frac{-1}{z^{1+q}}[\frac{\alpha_{0}}{1+q}+\frac{\alpha_{n}}{q}z+...+%
\frac{\alpha_{qn}}{1}z^{q}]$
\end{center}

\noindent avec ($\alpha_{o}$, $\alpha_{n}$,..., $\alpha_{qn})$ solution du
syst\`{e}me $(\mathcal{A)}$ et d\'{e}pendent seulement de $a_{n-1}$, $b_{m-s}
$ ($0$ $\leq$ $s$ $\leq$ $q$).\medskip

\underline{\textbf{2}}$^{eme}$\underline{\textbf{\ cas}}: $m=qn+r$, $q$ $\geq
$ $1$ et $o<r<n$.

\noindent Soit $J={\{}k$ tel que $0<k<m+n,k$ $\neq$ $m$ et $k$ $\notin$\ $n%
%TCIMACRO{\U{2115} }%
%BeginExpansion
\mathbb{N}
%EndExpansion
${\}}.

\noindent Si $k\in J$ alors $\alpha_{k}$ $=0$, en effet $\beta_{k,n}$ $=0$.

\noindent En comparant les termes de valuation comprise entre $-m-n$ et $%
-m-n+q+1$.

\noindent On a:

\begin{center}
$(\mathcal{B)}$ $\{%
\begin{array}{cc}
\beta_{sn,n}= & (-1)^{n}b_{m-s}\quad si\quad0\leq s\leq q+1 \\ 
\beta_{m,n}= & \alpha_{n-1}\beta_{0,n-1}%
\end{array}
$
\end{center}

\noindent Soit ($\alpha_{o}$, $\alpha_{n}$,..., $\alpha_{qn}$, $\alpha_{m}$, 
$\alpha_{(q+1)n})$ une solution du syst\`{e}me $(\mathcal{B)}$.

On a:

\begin{center}
$\alpha_{0}^{n}=(-1)^{n}b_{m}$ , $\alpha_{n}$ $=(-1)^{n}\frac{b_{m-1}}{%
n\alpha_{0}^{n-1}}$

$\alpha_{2n}=(-1)^{n}\frac{b_{m-2}}{n\alpha_{0}^{n-1}}-\frac{n-1}{2n^{2}}%
\frac{b_{m-1}^{2}}{\alpha_{0}^{2n-1}}$, ..., $\alpha_{m}=\frac{a_{n-1}}{n}%
.\medskip$
\end{center}

\textbf{1.2.3. }\underline{\textbf{Proposition}}: Si $m=qn+r$, $q$ $\geq$ $1 
$ et $0<r<n$ les facteurs d\'{e}terminants sont donn\'{e}s par:

\begin{center}
$Q(z)=$ $\frac{-1}{z^{1+\frac{m}{n}}}[\frac{n}{1+\frac{m}{n}}\alpha_{0}+%
\frac{n}{m}\alpha_{n}z+...+\frac{n}{n+r}\alpha_{qn}z^{q}+\alpha_{m}z^{\frac{m%
}{n}}+\frac{n}{r}\alpha_{(q+1)n}z^{q+1}]$
\end{center}

\noindent avec ($\alpha_{o}$, $\alpha_{n}$,..., $\alpha_{qn}$, $\alpha_{m}$, 
$\alpha_{(q+1)n})$ solution du syst\`{e}me $(\mathcal{B)}$ et d\'{e}pendent
seulement de $a_{n-1}$, $b_{m-s}$ ($0$ $\leq$ $s$ $\leq$ $q+1$).\medskip

\textbf{1.2.4}\textrm{\ }\textbf{P}\underline{\textbf{roposition}}: Soient $%
V_{1}$ et $V_{2}$ deux $D$-modules, d'Airy de bidegr\'{e} respectifs ($n_{1}$%
, $m_{1})$ et ($n_{2}$, $m_{2})$. On suppose que $V_{1}$ est formellement 
\'{e}quivalent \`{a} $V_{2}$ alors $n_{1}$ = $n_{2}$ et $m_{1}$= $m_{2}$%
.\medskip

\underline{\textbf{Preuve}}: Si on note par $\varphi$ l'isomorphisme
formelle entre $V_{1}$ et $V_{2}$ alors $\varphi$ transforme une base
formelle de $V_{1}$en une base formelle de $V_{2}$ on a alors $n_{1}$ = $%
n_{2}$. La pente \`{a} l'infini de $V_{1}$ est $\frac{n_{1}+m_{1}}{n_{1}}$
d'o\`{u} r\'{e}sultat.\medskip

Soient $A_{1}$ et $A_{2}$ deux op\'{e}rateurs d'Airy de bidegr\'{e} $(n,m)$, 
$V_{1}$ et $V_{2}$ leurs $D$-modules correspondant, on suppose que $V_{1}$
est formellement \'{e}quivalent \`{a} $V_{2}$. Les facteurs d\'{e}terminants
de $A_{1}$ et $A_{2}$ (en les num\'{e}rotants correctement) sont alors
identiques.

On d\'{e}signe par:

\begin{center}
$A_{1}=\sum_{i=1}^{n}a_{i}\partial^{i}-$ $\sum_{j=0}^{m}b_{j}x^{j}$ \quad($%
a_{n}$ $=1$)
\end{center}

et $A_{2}=\sum_{i=1}^{n}d_{i}\partial^{i}-\quad\sum_{j=0}^{m}h_{j}x^{j}$
\quad($d_{n}$ $=1$)\medskip

\textbf{1.2.5.} \underline{\textbf{Proposition}}: Si $\mathcal{D}$/$\mathcal{%
D}A_{1}$ est formellement \'{e}quivalent \`{a} $\mathcal{D}$/$\mathcal{D}%
A_{2}$ alors:

\noindent i) Si $m=qn$, $q$ $\geq$ $1$, on a:

\begin{center}
$a_{n-1}$ $=d_{n-1}$ et $b_{m-k}$ $=h_{n-k}$ pour $0$ $\leq$ $k$ $\leq$~$q$
\end{center}

\noindent ii) Si $m=qn+r$, $q$ $\geq$ 1 et $0<r<n$, on a:

\begin{center}
$a_{n-1}$ $=d_{n-1}$ et $b_{m-k}$ $=h_{m-k}$ pour $0$ $\leq$ $k$ $\leq$ $%
q+1.\medskip$
\end{center}

\underline{\textbf{Preuve}}: Les facteurs d\'{e}terminants sont invariant
par \'{e}quivalence formelle, les syst\`{e}mes $(\mathcal{A}$) et ($\mathcal{%
B}$) sont alors invariant, on utilise les propositions 1.2.2 et 1.2.3 pour
conclure.\medskip

\underline{\textbf{3}}$^{eme}$\underline{\textbf{\ cas}}: $n=qm+r$, $q$ $\geq
$ $1$ et $0<r<m$.

\noindent On a:

\[
\xi_{<o}=\frac{1}{z^{1+\frac{m}{n}}}[\alpha_{0}+\alpha_{m}z^{\frac{m}{n}%
}+...+\alpha_{(q+1)m}z^{(q+1)\frac{m}{n}}+\alpha_{n}z] 
\]

\noindent Dans $P_{L}$($z$,$\xi(z)$) $=0$, seuls les termes de valuation
comprise entre $-m-n$ et $-m-n+(q+1)m/n$ proviennent seulement des
coefficients d'un facteur d\'{e}terminant.\medskip

\textbf{1.2.6}\underline{\textbf{\ Proposition}}: $n=qm+r$, $q$ $\geq$ $1$, $%
0<r<$ $m$ et $\mathcal{D}$/$\mathcal{D}A_{1}$ formellement \'{e}quivalent 
\`{a} $\mathcal{D}$/$\mathcal{D}A_{1}$. Alors:

\noindent$b_{m}$ $=h_{m}$, $b_{m-1}$ $=h_{m-1}$

\noindent$a_{n-k}=d_{n-k}$ pour $0$ $\leq$ $k$ $\leq$ $q+1.\medskip$

\underline{\textbf{Preuve}}:

\[
\xi^{n+k}=\frac{1}{z^{n+m-k-k\frac{m}{n}}}[\beta_{n,n-k}z+(\sum_{j=0}^{q+1}%
\beta_{j,n-k}z^{j\frac{m}{n}})+\xi_{1} 
\]

\noindent o\`{u} $\xi_{1}$ est de valuation strictement sup\'{e}rieur \`{a} $%
(q+1)m/n$

\noindent on a:

$\beta_{o,n}$ $=(-1)^{n}$ $b_{m}$

$\beta_{n,n}$ $=(-1)^{n}$ $b_{m-1}$

\noindent et:

\[
(\mathcal{S)}\text{ }\mathcal{\{}%
\begin{array}{cc}
\beta_{1,n}-a_{n-1}\beta_{0,n-1}= & 0 \\ 
\beta_{2,n}-a_{n-1}\beta_{1,n-1}+a_{n-2}\beta_{0,n-2}= & 0 \\ 
\beta_{q+1,n}-a_{n-1}\beta_{q,n-1}+a_{n-2}%
\beta_{q-1,n-2}+..+(-1)^{q+1}a_{n-(q+1)}\beta_{0,n-(q+1)}= & 0%
\end{array}
\]

\noindent le syst\`{e}me $(\mathcal{S)}$ poss\'{e}de une solution unique ($%
a_{n-1}$,..., $a_{n-(q+1)})$ en fonction des $\beta_{i,j}$.\medskip

\underline{\textbf{Remarque}}: On obtient les r\'{e}sultats de cette
proposition par une transformation de Fourier de l'alg\`{e}bre de Weyl $%
%TCIMACRO{\U{2102} }%
%BeginExpansion
\mathbb{C}
%EndExpansion
\lbrack x$,$\partial]$.\medskip

\textbf{1.2.7. }\underline{\textbf{Monodromie formelle}}.\medskip

On pose $D=z\frac{d}{dz},\medskip$

on d\'{e}signe par $L=$ $\sum_{k=0}^{n}c_{k}$ $D^{n-k}$, $c_{o}$ $=1,$ la
forme de Frobenius- Fuchs de l'op\'{e}rateur d'Airy \S 1.2.

On pose pour $\alpha\in\bar{K},$ $L^{\alpha}=\sum_{k=0}^{n}c_{k}(D+%
\alpha)^{n-k}.$

Si $y$ est une solution de l$^{\text{'}}$\'{e}quation diff\'{e}rentielle $%
L(u)=0,$ alors $e^{-\int\frac{\alpha}{z}}y$ est une solution de l'\'{e}%
quation diff\'{e}rentielle $L^{\alpha}(v)=0.$

En effet:

\begin{center}
$(D+\alpha)^{k}(e^{-\int\frac{\alpha}{z}}y)=e^{-\int\frac{\alpha}{z}}D^{k}y.$
\end{center}

Si $Q$ est un facteur d\'{e}terminant pour $L$ donn\'{e} par:

\begin{center}
$z\frac{dQ}{dz}=\xi_{<0}$
\end{center}

\noindent on consid\`{e}re l'op\'{e}rateur diff\'{e}rentiel:

\begin{center}
$L^{\xi}=\sum_{k=0}^{n}c_{k}(D+\xi)^{n-k}=\sum_{k=0}^{n}h_{k}(\xi)D^{n-k}$
\end{center}

\noindent on a:

\begin{center}
$(D+\xi)^{n-k}=\sum_{i=0}^{n-k}\left( 
\begin{array}{c}
n-k \\ 
i%
\end{array}
\right) \xi^{\lbrack i]}$ $D^{n-k-i}$
\end{center}

\noindent o\`{u}:

\begin{center}
$\xi^{\lbrack k]}=e^{-\int\frac{\xi}{z}}$ $D^{k}(e^{\int\frac{\xi}{z}})$
\end{center}

\noindent on a:

\begin{center}
$\xi^{\lbrack k+1]}=\xi\xi^{\lbrack k]}$ $+D(\xi^{\lbrack k]})$
\end{center}

$\xi^{\lbrack0]}$ $=1$, $\xi^{\lbrack1]}=\xi$,\textrm{\ }$%
\xi^{\lbrack2]}=\xi^{2}+D\xi$

$\xi^{\lbrack3]}=\xi^{3}$ $+3\xi D\xi$ $+D^{2}\xi$, etc ....

\begin{center}
$L^{\xi}=$ $\sum_{k=0}^{n}c_{k}[\sum_{i=0}^{n-k}(%
\begin{array}{c}
n-k \\ 
i%
\end{array}
)\xi^{\lbrack i]}D^{n-k-i}]$
\end{center}

ce qui peut s'\'{e}crire sous la forme

\[
L^{\xi}=\sum_{k=0}^{n}[\sum_{i=0}^{k}c_{i}(%
\begin{array}{c}
n-i \\ 
k-i%
\end{array}
)\xi^{\lbrack k-i]}]D^{n-k} 
\]

\noindent et

\[
h_{k}(\xi)=\sum_{i=0}^{k}c_{i}(%
\begin{array}{c}
n-i \\ 
k-i%
\end{array}
)\xi^{\lbrack k-i]} 
\]

\noindent Il est facile de montrer, par r\'{e}currence, le r\'{e}sultat
suivant:

$\xi^{\lbrack k]}=\xi^{k}+\frac{k(k-1)}{2}\xi^{k-2}D\xi$ $+g_{k}(\xi$, $D\xi 
$,..., $D^{k-1}\xi)$

\noindent o\`{u} $g_{k}$ est un polyn\^{o}me en $k$ variables v\'{e}rifiant:

\begin{center}
$v(g_{k}(\xi$, $D\xi$,..., $D^{k-1}\xi)\geq$ $-(k-2)(1+$ $\frac{m}{n})$
\end{center}

\noindent on a:

\[
\{%
\begin{array}{c}
v(h_{k})=-k-k\frac{m}{n}\quad si\quad0\leq k\leq n-1 \\ 
v(h_{n})=v(h_{n-1})=-n-m+1+\frac{m}{n}%
\end{array}
\]

\noindent En effet:

\begin{center}
$h_{n}=$ $\sum_{k=0}^{n}c_{k}$ $\xi^{\lbrack n-k]}$
\end{center}

\noindent en utilisant l'\'{e}quation $P_{L}(z$, $\xi(z))=0$ et l'expression
des $\xi^{\lbrack k]}$, $k=0,...,$ $n.$

\noindent On obtient:

\[
\{%
\begin{array}{c}
h_{n}=\frac{\sigma}{z^{n+m-1-\frac{m}{n}}}+k_{n} \\ 
\sigma=\frac{1-n}{2}(n+m)\alpha_{0}^{n-1}\quad et\text{ }v(k_{n})>-n-m+\frac{%
m}{n}%
\end{array}
\]

\noindent D'autre part:

\begin{center}
$h_{n-1}=\sum_{i=0}^{n-1}c_{i}(%
\begin{array}{c}
n-i \\ 
n-1-i%
\end{array}
)\xi^{\lbrack n-1-i]}$
\end{center}

\noindent On a alors:

\[
\{%
\begin{array}{c}
h_{n-1}=\frac{n\alpha_{0}^{n-1}}{z^{n+m-1-\frac{m}{n}}}+l_{n}\text{ avec} \\ 
v(l_{n})>-n-m+1+\frac{m}{n}%
\end{array}
\]

\noindent$L^{\xi}$ est alors d'indice caract\'{e}ristique $n-1$, l'\'{e}%
quation indicielle est de degr\'{e} $1$.

Soit $w_{j}(z)=z^{-\lambda_{j}}v_{j}(z)$ une solution formelle de l'equation
diff\'{e}rentielle $L^{\xi}(u)=0$, on a:

\begin{center}
$-\lambda_{j}$ $n\alpha_{0}^{n-1}$ $+\sigma$ $=0$
\end{center}

\noindent d'o\`{u}

\begin{center}
$-\lambda_{j}$ $n\alpha_{0}^{n-1}+\frac{1-n}{2}(n+m)\alpha_{0}^{n-1}=0$
\end{center}

\noindent et

\begin{center}
$\lambda_{j}=\frac{1-n}{2n}$ $(n+m)$
\end{center}

\noindent$exp(2i\pi\lambda_{j})$ est une valeur propre de la partie de la
monodromie formelle attach\'{e}e \`{a} $Q$.

\begin{center}
$exp(2i\pi\lambda_{j})$ $=(-1)^{m+n-1}$ $exp(i\pi\frac{m}{n}).\medskip$
\end{center}

\noindent\textbf{2.}\underline{\textbf{Forme canonique d'une connexion d'Airy%
}}.\medskip

Soient $P_{n}$ et $Q_{m}$ deux polyn\^{o}mes de $%
%TCIMACRO{\U{2102} }%
%BeginExpansion
\mathbb{C}
%EndExpansion
\lbrack x]$ de degr\'{e} respectif $n$ et $m$ (on suppose que $P_{n}$ est
unitaire). On consid\`{e}re l'op\'{e}rateur diff\'{e}rentiel

\begin{center}
$L=P_{n}(\partial)$ $+Q_{m}(x)$
\end{center}

\noindent appell\'{e} op\'{e}rateur d'Airy de bidegr\'{e} $(n,m)$
correspondant \`{a} $P_{n}$ et $Q_{m}$ sur $\mathbb{P}^{1}(%
%TCIMACRO{\U{2102} }%
%BeginExpansion
\mathbb{C}
%EndExpansion
)$. \newline
Soit $\mathcal{F}$ le corps des s\'{e}ries de Laurent formelle, si $f\in%
\mathcal{F}$ alors $f=\sum_{j\in%
%TCIMACRO{\U{2124} }%
%BeginExpansion
\mathbb{Z}
%EndExpansion
}f_{j}z^{j}$ avec $f_{j}=0$ pour $\left\vert j\right\vert $ suffisament
grand. Soit $\overline{\mathcal{F}}$ la cl\^{o}ture alg\`{e}brique de $%
\mathcal{F}$, d'apr\`{e}s le th\'{e}or\`{e}me de Puiseux on a:

\[
\overline{\mathcal{F}}=\cup_{b\in%
%TCIMACRO{\U{2115} }%
%BeginExpansion
\mathbb{N}
%EndExpansion
_{\geq1}}\mathcal{F}(z^{\frac{1}{b}}) 
\]

\noindent L'op\'{e}rateur $L$ pr\'{e}sente une seule singularit\'{e}, elle
est situ\'{e}e \`{a} l'infini et elle est irr\'{e}guli\`{e}re. On pose $z=%
\frac {1}{x}$ et on se r\'{e}f\`{e}re \`{a}:

\[
A=\left( 
\begin{array}{ccccc}
0 & z^{-2} & 0 &  & 0 \\ 
0 & 0 & z^{-2} &  & 0 \\ 
&  &  &  &  \\ 
&  &  &  & z^{-2} \\ 
z^{-2}Q_{m}(\frac{1}{z}) & a_{1}z^{-2} &  &  & a_{n-1}z^{-2}%
\end{array}
\right) , 
\]

\noindent qu'on appelle connexion d'Airy associ\'{e}e \`{a} $(P_{n}$, $%
Q_{m}) $.

On se propose d'\'{e}tudier une classe d'\'{e}quivalence formelle de $A$,
c'est \`{a} dire l'ensemble des

$\xi\lbrack A]$ $=$ $\xi\mathrm{\ }A\xi^{-1}+\xi\prime\xi^{-1},$ o\`{u} $%
\xi\in Gl(n,\mathcal{\bar{F}})$ et $\xi\prime=\frac{d}{dz}\xi$.

On d\'{e}terminera une forme canonique de cette classe d'\'{e}quivalence et
les coefficients de $P_{n}$ et $Q_{m}$ qui interviennent dans une r\'{e}%
duction canonique de Babitt-Varadarajan [BV1].

Les facteurs d\'{e}terminants et la monodromie formelle de $L$ peuvent \^{e}%
tre calculer directement \`{a} partir de la forme canonique.

On d\'{e}finit ${\{}X_{n}$, $Y_{n}$, $H_{n}${\}} triplet standard par:

\[
X_{n}=\left( 
\begin{array}{ccccc}
0 & 1 & 0 &  & 0 \\ 
0 & 0 & 1 &  & 0 \\ 
&  &  &  &  \\ 
&  &  &  & 1 \\ 
0 &  &  &  & 0%
\end{array}
\right) 
\]

$Y_{n}=\left( 
\begin{array}{cccc}
0 &  &  & 0 \\ 
c_{1} &  &  &  \\ 
&  &  &  \\ 
&  & c_{n-1} & 0%
\end{array}
\right) $ (appell\'{e} nilpotent principal), les $c_{j}$ sont donn\'{e} par $%
c_{j}$ $=$ $j(n-j)$, $1$ $\leq$ $j$ $\leq$~$n-1.$

\noindent et $H_{n}=\left( 
\begin{array}{cccc}
n-1 &  &  &  \\ 
& n-3 &  &  \\ 
&  &  &  \\ 
&  &  & -(n-1)%
\end{array}
\right) $ les coefficients de la diagonale principale sont $%
h_{j,j}^{n}=n+1-2j.$

\noindent Si $H$ est une matrice semi-simple de $\mathcal{G}l(n,$ $%
%TCIMACRO{\U{2102} }%
%BeginExpansion
\mathbb{C}
%EndExpansion
),$n'ayant que des valeurs propres enti\`{e}res. On definit l'\'{e}l\'{e}%
ment $z^{H}$ de $Gl(n,$ $\mathcal{F})$ de la mani\`{e}re suivante: Si $%
\lambda\in%
%TCIMACRO{\U{2124} }%
%BeginExpansion
\mathbb{Z}
%EndExpansion
$ est une valeur propre de $H,$ et $V_{\lambda}$ est son sous espace propre
associ\'{e}, alors $z^{H}$ agit sur $V_{\lambda}$ par $z^{\lambda}I.$ L'\'{e}%
lement $z^{H}$ commute avec tous les endomorphismes de $%
%TCIMACRO{\U{2102} }%
%BeginExpansion
\mathbb{C}
%EndExpansion
^{n}$ qui commutent avec $H.$ Si les valeurs propres de $H$ sont des \'{e}l%
\'{e}ments de $(\frac{1}{q})%
%TCIMACRO{\U{2124} }%
%BeginExpansion
\mathbb{Z}
%EndExpansion
,$ avec $q\in%
%TCIMACRO{\U{2115} }%
%BeginExpansion
\mathbb{N}
%EndExpansion
_{\geq1}.$ En \'{e}crivant $z^{\frac{\mu}{q}}=(z^{\frac{1}{q}})^{\mu},$ $%
\mu\in%
%TCIMACRO{\U{2124} }%
%BeginExpansion
\mathbb{Z}
%EndExpansion
,$ on d\'{e}finit de la m\^{e}me mani\`{e}re $z^{H}.$ On a, par exemple: $%
z^{H_{n}}=\left( 
\begin{array}{cccc}
z^{n-1} & 0 &  & 0 \\ 
& z^{n-3} &  &  \\ 
0 &  &  &  \\ 
0 & 0 &  & z^{-(n-1)}%
\end{array}
\right) .$

\noindent En \'{e}crivant la matrice $A$ sous la forme: $A=%
\sum_{k=-m-2}^{-2}A_{k}z^{k},$ on remarque que $A_{-m-2}=\left( 
\begin{array}{cccc}
0 & 0 &  & 0 \\ 
&  &  &  \\ 
0 &  &  & 0 \\ 
b_{m} & 0 &  & 0%
\end{array}
\right) $ qui est une matrice nilpotente. Le niveau principal (o\`{u}
invariant principal), voir [BV1], est donc strictement sup\'{e}rieur \`{a} $%
-m-2.$ La transformation de "choix" (ou "Shearing") $z^{\frac{s}{2}H_{n}}$
appliqu\'{e}e \`{a} $A$ donne un \'{e}l\'{e}ment de la classe d'\'{e}%
quivalence formelle:

\[
z^{\frac{s}{2}H_{n}}[A]=\left( 
\begin{array}{cccc}
(n-1)\frac{s}{2}z^{-1} & z^{-2+s} &  & 0 \\ 
& (n-3)\frac{s}{2}z^{-1} &  &  \\ 
0 &  &  & z^{-2+s} \\ 
z^{-2-(n-1)s}Q_{m}(\frac{1}{z}) & a_{1}z^{-2-(n-2)s} &  & a_{n-1}z^{-2}-(n-1)%
\frac{s}{2}z^{-1}%
\end{array}
\right) 
\]

\noindent Le choix $s=-\frac{m}{n},$ donne:

\begin{center}
$(\ast)$ $z^{\frac{s}{2}H_{n}}[A]=A^{1}=$ $z^{-\frac{m}{n}-2}A_{-\frac{m}{n}%
-2}^{1}+$ $z^{-\frac{m}{n}-1}A_{-\frac{m}{n}-1}^{1}+$ $...$
\end{center}

\noindent avec $A_{-\frac{m}{n}-2}^{1}=\left( 
\begin{array}{cccc}
0 & 1 &  & 0 \\ 
&  &  &  \\ 
0 &  & 0 & 1 \\ 
b_{m} & 0 &  & 0%
\end{array}
\right) $

\noindent qui n'est pas nilpotente, on a alors l'invariant principal de $A$
est $r=-\frac{m}{n}-2.$

\noindent L'op\'{e}rateur $A_{-\frac{m}{n}-2}^{1}$ est semi-simple et admet $%
n$ valeurs propres distinctes $e^{\frac{2i\pi k}{n}}b_{m}^{\frac{1}{n}},$ $%
k=1,...,n$ (moyennant un choix de $b_{m}^{\frac{1}{n}}).$

On note $\mathcal{G=G}l(n,$ $%
%TCIMACRO{\U{2102} }%
%BeginExpansion
\mathbb{C}
%EndExpansion
)$ l'alg\`{e}bre de Lie de $GL(n,%
%TCIMACRO{\U{2102} }%
%BeginExpansion
\mathbb{C}
%EndExpansion
)$. On consid\`{e}re l'action adjointe de $GL(n,%
%TCIMACRO{\U{2102} }%
%BeginExpansion
\mathbb{C}
%EndExpansion
)$ sur $\mathcal{G}$, on a:

\begin{center}
$\mathcal{G=G}_{A_{r}^{1}}\oplus\lbrack A_{r}^{1},\mathcal{G]}$
\end{center}

\noindent o\`{u} $\mathcal{G}_{A_{r}^{1}}$ est le commutateur de $A_{r}^{1}. 
$

\noindent$ad$ $A_{r}^{1}$ est un isomorphisme sur $[A_{r}^{1},\mathcal{G]}$;
utilisant les notations de Babitt-Varadarajan [BV1] on a:

\noindent si $V$ est un espace vectoriel complexe de dimension $n$

$V(\mathcal{F})=$ $V\otimes_{%
%TCIMACRO{\U{2102} }%
%BeginExpansion
\mathbb{C}
%EndExpansion
}\mathcal{F},$ $\mathcal{G}l(V(\mathcal{F}))=\mathcal{G}l(V)\otimes_{%
%TCIMACRO{\U{2102} }%
%BeginExpansion
\mathbb{C}
%EndExpansion
}\mathcal{F}$ et $GL(V(\mathcal{F}))$ $=Aut_{\mathcal{F}}$($V(\mathcal{F})$).

\noindent On note par $\mathcal{O}$ l'anneau des s\'{e}ries de Laurent
formelle d'ordre $\geq$ $0$. Si $p\in%
%TCIMACRO{\U{2115} }%
%BeginExpansion
\mathbb{N}
%EndExpansion
_{\geq1},$ $GL(V(\mathcal{O}))_{p}$ d\'{e}signe le sous groupe des \'{e}l%
\'{e}ment de $GL(V$($\mathcal{O}))$ congru a $I$ modulo $z^{p}$.\medskip

\textbf{2.1. }\underline{\textbf{Lemme [BV1]}}: (projection spectrale)

Soit $A=A_{r}$ $z^{r}$ $+$ $A_{r+1}$ $z^{r+1}$ $+$... , $r<-1,$un \'{e}l\'{e}%
ment de $\mathcal{G}l(V(\mathcal{F}))$, $\Sigma$ l'ensemble des valeurs
propres de $A_{r}$; pour $\lambda\in\Sigma$, $P_{\lambda}$ la projection
spectrale correspondante. Il existe alors\textrm{\ }$\Psi\in GL(V$($\mathcal{%
O}))_{1}$ tel que $A^{(1)}$ $=$ $\Psi\lbrack A]$ commute avec tous les $%
P_{\lambda}$; en effet si $X$ est une composante semi-simple de $A_{r}$, il
existe $T_{k}$ unique dans $[X,\mathcal{G}]$ ($k$ $\geq$ $1$) tel que si:

\begin{equation}
\Psi=\sqcap_{k=1}^{\infty}(I+z^{k}T_{k})=Lim_{m\rightarrow%
\infty}(I+z^{m}T_{m})\text{ ... }(I+zT_{1})   \label{eq1}
\end{equation}

\noindent alors $A^{(1)}$ = $\Psi\lbrack A]$ commute avec tous les $%
P_{\lambda}\mathrm{.}$

Si $Y\in$ $GL(V(\mathcal{O}))_{1}$ est tel que $A^{(2)}=Y[A]$ commute avec
tous les $P_{\lambda}$ alors $Y=U\Psi$ et $A^{(2)}$ $=U[A^{(1)}]$ o\`{u} $%
U\in$ $GL(V(\mathcal{O}))_{1}$ commute avec tous les $P_{\lambda }\mathrm{%
.\medskip}$

Soit $A$ une connexion de $\mathcal{G}l(V(\mathcal{\bar{F}}))$ d\'{e}finie
par: $A=\sum_{k=s}^{\infty}A_{r_{k}}z^{r_{k}},$ avec $r_{k}\in%
%TCIMACRO{\U{211a} }%
%BeginExpansion
\mathbb{Q}
%EndExpansion
,$ $A_{r_{k}}$ $\neq0$ et $r_{s}<-1.$ On appelle niveau de $A$ la suite $%
(r_{s}$, $r_{s+1}$,..., $r_{m})$ o\`{u} $m$ est l'unique entier v\'{e}%
rifiant $r_{m+1}\quad\geq$ $-1$.\medskip

\underline{\textbf{Exemple}}: Si $m$ et $n$ sont deux entiers naturels non
nul et $A=A_{-\frac{m}{n}-2}z^{-\frac{m}{n}-2}+A_{-2}z^{-2}+A_{-1}z^{-1}+...%
\quad,$ avec $A_{-\frac{m}{n}-2}$ et $A_{-2}$ non nuls alors le niveau de $A 
$ est par d\'{e}finition le couple $(-\frac{m}{n}-2,$ $-2).\medskip$

\noindent En g\'{e}n\'{e}ral, les transformations successives $1+$ $%
z^{k}T_{k}$ donn\'{e}es dans le lemme peuvent changer le niveau d'une
connexion.\medskip

\textbf{2.2. }\underline{\textbf{D\'{e}finition}}: On appelle mod\`{e}le
canonique une connexion de la forme:

$B=D_{r_{1}}z^{r_{1}}+D_{r_{2}}z^{r_{2}}+...$ $+D_{r_{m}}z^{r_{m}}+$ $%
Cz^{-1},$ avec $r_{1}<$ $r_{2}$ $<$...$<$ $r_{m}$ $<$ $-1.$ Les op\'{e}%
rateurs $D_{r_{1}},$ ..., $D_{r_{m}}$ et $C$ sont des \'{e}l\'{e}ments de $%
\mathcal{G}l(V)$ qui commutent deux \`{a} deux et les $D_{r_{j}}$ sont
semi-simples.\medskip

Babitt-Varadajan [BV1] ont montr\'{e} que toute connexion dans $\mathcal{G}%
l(V(\mathcal{\bar{F}}))$ est \'{e}quivalente \`{a} un mod\`{e}le canonique.
Le niveau ($r_{1}$,..., $r_{m})$ du mod\`{e}le canonique est invariant par 
\'{e}quivalence formelle, en particulier $r_{1}$ est l'invariant principal.
Deux mod\`{e}les canoniques d\'{e}finis par:

\[
(D_{r_{1}},\text{..., }D_{r_{m}},\text{ }C)\text{ et }(D_{r_{1}}^{\prime },%
\text{..., }D_{r_{m}}^{\prime},\text{ }C^{\prime}) 
\]

\noindent sont dans la m\^{e}me classe d'\'{e}quivalence formelle si et
seulement si pour $k$ $\geq$ $1$

\[
(D_{r_{1}},\text{..., }D_{r_{m}},\exp2\pi ik\text{ }C)\text{ et }%
(D_{r_{1}}^{\prime},\text{..., }D_{r_{m}}^{\prime},\exp2\pi ik\text{ }%
C^{\prime}) 
\]

\noindent sont dans la m\^{e}me $GL(V)$-orbite.\medskip

Dans la suite on \'{e}tudiera un mod\`{e}le canonique d'une connexion
d'Airy, on d\'{e}terminera les coefficients $a_{i}$ et $b_{j}$ de l'op\'{e}%
rateur qui interviennent dans la r\'{e}duction canonique.\medskip

Soit $A^{1}$ la connexion de $\mathcal{G}l(V(\mathcal{\bar{F}}))$ donn\'{e}e
dans \S 2. $(\ast)$ par:

$A^{1}=A_{r}^{1}$ $z^{r}$ + $A_{r+1}^{1}z^{r+1}+...$ , avec $r=-\frac{m}{n}%
-2.$ L'espace $\mathcal{G}_{A_{r}^{1}}$ est d\'{e}finie par:

\[
K\in\mathcal{G}_{A_{r}^{1}}\Longleftrightarrow K=\left( 
\begin{array}{cccc}
\alpha_{1} & \alpha_{2} &  & \alpha_{n} \\ 
b_{m}\alpha_{n} & \alpha_{1} &  &  \\ 
&  &  & \alpha_{2} \\ 
b_{m}\alpha_{2} &  & b_{m}\alpha_{n} & \alpha_{1}%
\end{array}
\right) ,\text{ avec }\alpha_{i}\in%
%TCIMACRO{\U{2102} }%
%BeginExpansion
\mathbb{C}
%EndExpansion
. 
\]

\noindent$ad_{A_{r}^{1}}$ est inversible sur $[A_{r}^{1}$, $\mathcal{G]},$
il existe alors un op\'{e}rateur $T_{1}$ unique dans $[A_{r}^{1}$, $\mathcal{%
G]}$ tel que:

\begin{center}
($I+zT_{1})$ $[A^{1}]$ $=A^{2}$ $=A_{r}^{1}$ $z^{r}$ + $%
A_{r+1}^{2}z^{r+1}+...,$ avec

$A_{r+1}^{2}=A_{r+1}^{1}-[A_{r}^{1},$ $T_{1}]$ $\in\mathcal{G}_{A_{r}^{1}}$
\end{center}

\noindent en effet, $A^{1}=A_{r}^{1}$ $z^{r}$ + $A_{r+1}^{1}z^{r+1}+...$ , l'%
\'{e}quation ($I+zT_{1})$ $[A^{1}]$ $=A^{2}$ est \'{e}quivalente \`{a}:

\begin{center}
($I+zT_{1})(A_{r}^{1}$ $z^{r}$ + $A_{r+1}^{1}z^{r+1}+...)$ $=A^{2}$($%
I+zT_{1})$ $-T_{1}$
\end{center}

\noindent en \'{e}crivant:

\begin{center}
$A^{2}=A_{r}^{2}$ $z^{r}$ + $A_{r+1}^{2}z^{r+1}+...$
\end{center}

\noindent et en comparant les coefficients des $z^{i},$ $i\geq r,$ on
obtient: $A_{r}^{1}=A_{r}^{2}$ et $A_{r+1}^{2}=A_{r+1}^{1}$ $-$ $[A_{r}^{1},$
$T_{1}].$ Comme $\mathcal{G=G}_{A_{r}^{1}}\oplus\lbrack A_{r}^{1},\mathcal{G]%
},$ et $ad_{A_{r}^{1}}$ est inversible sur $[A_{r}^{1},\mathcal{G]}$, il
existe un unique $T_{1}$ dans $[A_{r}^{1},\mathcal{G]}$ tel que $%
A_{r+1}^{2}\in\mathcal{G}_{A_{r}^{1}}.$

\noindent on remarquera que si $b_{m-1}\neq$~$0$ alors

\begin{center}
$A_{r+1}^{2}\neq0$ car sinon $A_{r+1}^{1}=\left( 
\begin{array}{cccc}
0 & 0 & 0 & 0 \\ 
0 & 0 & 0 & 0 \\ 
0 & 0 & 0 & 0 \\ 
b_{m-1} & 0 & 0 & 0%
\end{array}
\right) $
\end{center}

\noindent serait dans $[A_{r}^{1}$, $\mathcal{G]}$ ce qui est absurde. Ceci
nous permet d\'{e}j\`{a} de remarquer que si $r+1$ figure dans le niveau de $%
A^{1}$, il figure aussi dans le niveau de $A^{2}=(1+zT_{1})$ $[A^{1}]$. On
peut calculer exactement $T_{1}$ et donner la forme explicite de $A^{2}$.

\noindent L'op\'{e}rateur $T_{1}$ est donn\'{e}e par:

\[
T_{1}=\frac{b_{m-1}}{nb_{m}}\left( 
\begin{array}{cccc}
0 & 0 &  & 0 \\ 
& -1 &  & 0 \\ 
0 &  &  &  \\ 
0 & 0 &  & -(n-1)%
\end{array}
\right) , 
\]

\[
A_{r+1}^{2}=\frac{b_{m-1}}{nb_{m}}\left( 
\begin{array}{cccc}
0 & 1 &  & 0 \\ 
& 0 & 1 &  \\ 
0 &  &  & 1 \\ 
b_{m} &  &  & 0%
\end{array}
\right) . 
\]

\noindent En g\'{e}n\'{e}ral, pour $k\geq1,$ on a:

\[
A_{r+k}^{2}=A_{r+k}^{1}+T_{1}A_{r+k-1}^{1}-A_{r+k-1}^{2}T_{1}. 
\]

\noindent On peut, aussi, retrouver directement les coefficients de $A^{2}$%
par:

\begin{center}
$(1+zT_{1})$ $[A^{1}]=A^{2}$

\[
A^{2}=\left( 
\begin{array}{cccc}
0 & p_{1}(z) & 0 & 0 \\ 
&  &  & 0 \\ 
0 &  & 0 & p_{n-1}(z) \\ 
q_{0}(z) & q_{1}(z) &  & q_{n-1}(z)%
\end{array}
\right) -\frac{m}{2n}z^{-1}H_{n}+.. 
\]
\end{center}

\noindent avec:

$p_{j}(z)=z^{-\frac{m}{n}-2}(1-(j-1)z\frac{b_{m-1}}{nb_{m}})(1-jz\frac {%
b_{m-1}}{nb_{m}})^{-1},$ $1\leq j\leq n-1,$

$q_{0}(z)=z^{-2+\frac{m}{n}(n-1)}Q_{m}(\frac{1}{z})(1-(n-1)z\frac{b_{m-1}}{%
nb_{m}}),$

$q_{k}(z)=a_{k}z^{-2+\frac{m}{n}(n-1-k)}(1-(n-1)z\frac{b_{m-1}}{nb_{m}})(1-kz%
\frac{b_{m-1}}{nb_{m}})^{-1},$ $1\leq k\leq n-1,$

$q_{n-1}(z)=a_{n-1}z^{-2}$ \medskip

\noindent En d\'{e}veloppant au voisinage de $0$ les fonctions $p_{i}$ et $%
q_{j}$ on obtient les coefficients de $A^{2}$.

Pour construire d'autre transformation du type $1+z^{k}$ $T_{k}$, $k$ $\geq$ 
$1$, il est n\'{e}cessaire de consid\'{e}rer les cas $m>n$, $m<n$ et $m=n$. 
\newline
L'\'{e}tude des coefficients de la connexion $A$ montre que les coefficients
qui peuvent intervenir dans le mod\`{e}le canonique dans le cas $n=qm+s$, $q$
$\geq$ $1$, $0$ $<s<m$; sont $b_{m}$, $b_{m-1}$, $a_{n-1}$,..., $a_{n-(q+1)}$%
. Dans le cas $m=qn+s$, $q$ $\geq$ $1$, $0$ $<s<m$, sont $a_{n-1}$, $b_{m}$%
,..., $b_{m-(q+1)}$. Enfin, dans le cas $m=qn$ sont $b_{m}$, $b_{m-1}$, ..., 
$b_{m-q}$, $a_{n-1}$.\medskip

\textbf{2.3. }\underline{\textbf{Proposition}} :

Soit $A$ une connexion d'Airy de bidegr\'{e} $(n,m)$, $m=nq+s$, $q$ $\geq$ $1
$, $0$ $<s<n$. Il existe alors des op\'{e}rateurs $T_{k}$, $1$ $\leq$ $k$ $%
\leq$~$q+1$, et $T_{\frac{m}{n}}$ unique dans $[A_{r},\mathcal{G}]$, tel que
la connexion:

\[
A^{q+1}=(I+z^{q+1}T_{q+1})(I+z^{\frac{m}{n}}T_{\frac{m}{n}})(I+z^{q}T_{q})%
\text{ }...(I+zT_{1})[A] 
\]

\noindent soit donn\'{e}e par :

\begin{center}
$A^{q+1}$ $=$ $(z^{r}+\alpha_{1}$ $z^{r+1}$ $+$...$+$ $\alpha_{q+1}$ $%
z^{r+q+1})$ $A_{r}+$ $z^{-2}\frac{a_{n-1}}{n}I$ + $z^{-1}\left( \ast\right) $
+...
\end{center}

\noindent Les coefficients $\alpha_{1},...,\alpha_{q+1}$ d\'{e}pendent
seulement de $b_{m}$,..., $b_{m-(q+1)}.\medskip$

\underline{\textbf{Preuve}}: On \'{e}tablit, d'abord, l'existence de $T_{1}$%
. Puis, de la m\^{e}me mani\`{e}re on montre l'existence des $T_{k}$, ($2$ $%
\leq$ $k$ $\leq$ $q$), les transfrormations $I$+$z^{k}T_{k}$ conservent les
coefficients de $z^{r+t}$ ( $0$ $\leq$ $t$ $\leq$ $k-1$) et changent le
coefficient de $z^{r+k}$ en un \'{e}l\'{e}ment de $\mathcal{G}_{A_{r}}$ et
modifient les coefficients de $z^{r+t}$ ( $t$ $\geq$ $k+1$). Enfin les op%
\'{e}rateurs $T_{k}$ sont diagonaux et d\'{e}pendent seulement de $b_{m}$%
,..., $b_{m-k}$.

\noindent En effet le coefficient de $z^{r+k}$ dans $A^{k-1}$ est de la forme

\[
A_{r+k}^{k-1}=\left( 
\begin{array}{cccc}
0 & \ast & 0 & 0 \\ 
& 0 & \ast & 0 \\ 
0 & 0 &  & \ast \\ 
\ast & 0 &  & 0%
\end{array}
\right) 
\]

\noindent seul les coefficients de la $i$-\`{e}me ligne ($i+1)$-\`{e}me
colonne (($i,i+1))$, et $n-$\`{e}me ligne premi\`{e}re colonne $(n,1)$ sont
eventuellement non nul. On a,

\begin{center}
$A_{r+k}^{k}=A_{r+k}^{k-1}-[A_{r}$, $T_{k}]$
\end{center}

\noindent il existe $T_{k}$ unique dans $[A_{r},$ $\mathcal{G}$], tel que
l'op\'{e}rateur $A_{r+k}^{k}\in\mathcal{G}_{A_{r}}.$ L'op\'{e}rateur $T_{k}$
est diagonale (seul les \'{e}l\'{e}ments de la diagonale principale
interviennent), et une it\'{e}ration simple montre que $T_{k}$ d\'{e}pend
seulement de $b_{m}$, .., $b_{m-k}$. L'hypoth\'{e}se $m=nq+s,$ $s$ $\neq$ $0 
$, entraine directement que les coefficients de $z^{-1}$ et $z^{-2}$ ne sont
pas affect\'{e}s par la transformation $I$+$z^{k}T_{k}$. \newline
On peut toujours choisir une branche de $\xi^{n}$ $=z$, et choisir $X$ dans $%
GL(V(\mathcal{O}))_{1}$ de fa\c{c}on a r\'{e}duire la connexion $A$.

Soit

\begin{center}
$A^{q}$ $=(I+z^{q}T_{q})$...$(I+zT_{1})$ $[A]$

$=(z^{r}+\alpha_{1}z^{r+1}$+...+ $\alpha_{q}z^{r+q})A_{r}$ + $z^{-2}\left( 
\begin{array}{ccc}
0 &  & 0 \\ 
0 &  & 0 \\ 
0 &  & a_{n-1}%
\end{array}
\right) +z^{-1-\frac{s}{n}}\left( 
\begin{array}{ccc}
0 & \ast &  \\ 
&  & \ast \\ 
\ast &  & 0%
\end{array}
\right) -\frac{m}{2n}z^{-1}H_{n}$+...
\end{center}

\noindent Il existe, dans $[A_{r}$, $\mathcal{G}$], un op\'{e}rateur $T_{%
\frac{m}{n}}$ unique tel que le coefficient $A_{-2}^{\frac{m}{n}}$ de $z^{-2}
$ dans ($I$ + $z^{\frac{m}{n}})[A^{q}]$ = $A^{\frac{m}{n}}$ soit dans $%
\mathcal{G}_{A_{r}}$,

\begin{center}
$A_{-2}^{\frac{m}{n}}=[$ $T_{\frac{m}{n}},A_{r}]+\left( 
\begin{array}{ccc}
0 & ... & 0 \\ 
&  &  \\ 
0 & ...0 & a_{n-1}%
\end{array}
\right) \in\mathcal{G}_{A_{r}}$
\end{center}

\noindent et $A_{-2}^{\frac{m}{n}}=\frac{a_{n-1}}{n}I.$

\noindent Le coefficient de $z^{-1-\frac{s}{n}}$ dans $A^{\frac{m}{n}},$
soit $A_{-1-\frac{s}{n}}^{\frac{m}{n}}$ = $A_{-1-\frac{s}{n}}^{q}$et $%
A_{-1}^{\frac{m}{n}}$ $=\frac{b_{m-1}}{nb_{m}}[T_{\frac{m}{n}},A_{r}]$ $-%
\frac {m}{2n}H_{n}$ est une matrice diagonale.

\noindent Il existe $T_{q+1}$ unique dans $[A_{r}$, $\mathcal{G}$] tel que
le coefficient $A_{-1-\frac{s}{n}}^{q+1}$ de $z^{-1-\frac{s}{n}}$ dans ($I$+$%
z^{q+1}$ $T_{q+1})$ $[A^{\frac{m}{n}}]$ = $A^{q+1}$ est dans $\mathcal{G}%
_{A_{r}}$.

\begin{center}
$A_{-1-\frac{s}{n}}^{q+1}=A_{-1-\frac{s}{n}}^{\frac{m}{n}}+[T_{q+1},A_{r}],$
avec

$A_{-1-\frac{s}{n}}^{\frac{m}{n}}=\left( 
\begin{array}{ccc}
0 & \ast & 0 \\ 
&  & \ast \\ 
\ast & 0.. & 0%
\end{array}
\right) .$
\end{center}

\noindent L'op\'{e}rateur $T_{q+1}$ est alors diagonal et $A_{-1-\frac{s}{n}%
}^{q+1}=\alpha_{q+1}$ $A_{r}$.\medskip

\textbf{2.4. }\underline{\textbf{Th\'{e}or\`{e}me}}: Il existe $U$ dans $%
GL(V(\mathcal{\bar{F}}))$ tel que:

\begin{center}
$U[A]$ $=B=diag(q_{1}(z)$,...,$q_{n}(z)$) + $\frac{1-n}{2n}(n+m)z^{-1}I,$
\end{center}

\noindent o\`{u} $q_{1}$,...,$q_{n}$ sont des polyn\^{o}mes en $z^{\frac {-1%
}{n}}$ de degr\'{e} $m+2n$ et dont les coefficients sont des fonctions de $%
a_{n-1}$, $b_{m}$, $b_{m-1}$,..., $b_{m-(q+1)}.\medskip$

\underline{\textbf{Preuve}} :

On applique la proposition 2.3 puis le lemme 2 [BVI] page 44 et le th\'{e}or%
\`{e}me [BV1] page 45, la connexion $A$ est alors \'{e}quivalente sur $%
\mathcal{\bar{F}}$ \`{a}:

\begin{center}
$B^{1}$ $=(z^{r}+\alpha_{1}z^{r+1}$ +...+ $\alpha_{q+1}$ $z^{r+q})A_{r}+%
\frac{a_{n-1}}{n}z^{-2}I$ + $z^{-1}C$
\end{center}

\noindent L'op\'{e}rateur $A_{r}$ est semi simple et commute avec $C$, du
pragraphe 1.2.7, on d\'{e}duit que $\frac{1-n}{2n}(n+m)$ est l'unique valeur
propre de $C$. Il existe alors une base de $V$ dans laquelle $A_{r}$ est
diagonal et $C$ est \'{e}quivalente \`{a} $(\frac{1-n}{2n})(n+m)I$.\medskip

\underline{\textbf{Remarque}}: Les facteurs d\'{e}terminants sont donn\'{e}%
es par:

$Q_{i}(z)$ = $\int q_{i}(z)dz$ (la constante d'int\'{e}gration est nul).

La monodromie formelle est donn\'{e}e par $exp(2i\pi C)$.\bigskip

\begin{center}
\textbf{R\'{e}f\'{e}rences\medskip}
\end{center}

\noindent\lbrack BV1] Babitt, (D.G); Varadarajan, (V.S): "Formal Reduction
of Meromorphic Differential: A Group theoretic View", Pacific J. Math, 108,
1-80, (1983).\medskip

\noindent\medskip\lbrack BV2] Babitt, (D.G); Varadarajan, (V.S):
"Deformation Of Nilpotent Matrices Over Rings And reduction Of Analytic
Families Of Meromorphic Differential Equations", Mem. Amer. Soc., 55, N$^{0}$
325, (1985).\medskip

\noindent\lbrack F] Forsyth, (A.\ R): "Theory of differential Equations",
Cambridge Univ. Press. London, New-York, (1890).\medskip

\noindent\lbrack H K] Helfer, (B.); Kannai, (Y.): "Determining Factors and
Hypoellipticity of Ordinary Differential Operators With Double
'characteristics'", Ast\'{e}risque 2.3, 198-216,(1973).\medskip

\noindent\lbrack S1] Sa\"{\i}dane, (L.): "Propri\'{e}t\'{e}s Alg\'{e}briques
des Op\'{e}rateurs d'Airy de petit Ordre", Ann. Fac. Sc. Toulouse, Vol IX, n$%
^{0}$ 3, 519-550, (2000).\medskip

\noindent\lbrack S2] Sa\"{\i}dane, (L.): "Crit\`{e}res de r\'{e}ductibilit%
\'{e} et d'\'{e}quivalence pour les op\'{e}rateurs d'Airy de petit \qquad
ordre", C. R. Acad. Sci. Paris, t. 321, S\'{e}rie I, 523 - 526, 1995.\medskip

\noindent \lbrack Y] Yebbou, (J.): "Calcul des facteurs d\'{e}terminants",
J. Diff. Equations 72, 140-148, (1988).

\end{document}